\documentclass[11pt]{article}
\usepackage{amsfonts,graphics}

\usepackage{cite}
\usepackage{latexsym}
\usepackage{amsmath}
\usepackage{amssymb}
\usepackage[dvips]{epsfig}
\usepackage{amscd}

\newtheorem{theorem}{Theorem}[section]
\newtheorem{remark}{Remark}[section]
\newcommand{\thatsall}{\hfill$\Box$}
\newtheorem{lemma}{Lemma}[section]

\newcommand{\ti}{\tilde}
\def\pf{{\it Proof.}  }

\newcommand{\te}{\theta}

\newcommand{\pa}{\partial}
\def\O{\Omega}
\def\r{\mathbb{R}}
\newcommand{\bi}{\bibitem}

\newcommand{\dis}{\displaystyle}

\renewcommand{\b}{\beta}
\newcommand{\bl}{\begin{lemma}}
	\newcommand{\el}{\end{lemma}}
\newcommand{\et}{\end{theorem}}
\newcommand{\ga}{\gamma}

\def\xix{\int_{\O} }

\newcommand{\al}{\alpha}
\newcommand{\de}{\delta}
\newcommand{\ve}{\varepsilon}
\newcommand{\la}{\label}
\newcommand{\ka}{\kappa}

\newcommand{\ol}{\overline}

\newcommand{\bn}{\begin{eqnarray}}
\newcommand{\en}{\end{eqnarray}}
\newcommand{\bess}{\begin{eqnarray*}}
\newcommand{\eess}{\end{eqnarray*}}

\newcommand{\be}{\begin{equation}}
\newcommand{\ee}{\end{equation}}

\newcommand{\ba}{\begin{aligned}}
\newcommand{\ea}{\end{aligned}}
\newcommand{\bnn}{\begin{eqnarray*}}
\newcommand{\enn}{\end{eqnarray*}}

\def\norm[#1]#2{\|#2\|_{#1}}

\def\a{\alpha}

\def\O{\Omega}
\def\rr{\mathbb{R}}

\renewcommand{\b}{\beta  }

\makeatletter      
\@addtoreset{equation}{section}
\makeatother       
\date{}

\usepackage{url}
\newcommand{\authoraddress}[2]{%
\textsc{#1}, \textit{E-mail address:}  \protect\url{#2}
}
\newcommand{\KLAddress}{
\authoraddress{Department of Math, University of Michigan, United States}{kexinli@umich.edu}%
}
\newcommand{\XXAddress}{%
\authoraddress{School of Mathematical Sciences, Beijing Normal University, China}{xjxu@bnu.edu.cn}, \textit{Funding:}  Xu was partially supported by the National Key R\&D Program of China (grant 2020YFA0712900) and the National Natural Science Foundation of China (grants 12171040, 11771045, and 12071069)  %
}

\title{Large-Time Behavior of Solutions to
Compressible Navier-Stokes
System in Unbounded Domains with Degenerate Heat-Conductivity  and Large Data}
\author{Kexin L{\small I}\thanks{\KLAddress},\quad Xiaojing X{\small U}\thanks{\XXAddress}   }

\date{}

\begin{document}
\maketitle



\begin{abstract} We are concerned with the large-time behavior of solutions
	to the initial and initial boundary value problems with  large initial data for the compressible Navier-Stokes system with degenerate heat-conductivity  describing  the  one-dimensional motion of a
	viscous heat-conducting perfect polytropic  gas in unbounded domains.
	Both the specific volume and   temperature are  proved to be  bounded  from below and above independently of both time and space. Moreover, it is shown that the global solution is    asymptotically stable as time tends to infinity.  \end{abstract}


\section{Introduction}

Describing the  one-dimensional motion of a
viscous heat-conducting  polytropic  gas, the compressible Navier-Stokes system
is written in the Lagrange variables in the following form (see \cite{ba,se})
\be \la{1.1}
v_t=u_{x},
\ee
\be \la{1.2}
u_{t}+P_{x}=\left(\mu\frac{u_{x}}{v}\right)_{x},\ee
\be  \la{1.'3} \left(e+\frac{u^2}{2}\right)_{t}+ (P
u)_{x}=\left(\kappa\frac{\theta_{x}}{v}+\mu\frac{uu_x}{v}\right)_{x}
,
\ee   where   $t>0$ is time,  $x\in\O\subset\rr=(-\infty,\infty)$ denotes the
Lagrange mass coordinate,   the unknown functions $v ,u,$ $\theta ,e ,$ and $P$ are,  respectively, the specific volume of the gas,  fluid velocity,  absolute temperature, internal energy, and  pressure;
$\mu $ is the viscosity coefficient, and  $\ka $  is the
heat conductivity one. In general, $P,e,$    $\mu$, and  $\ka$ are functions of $\te$ and
$v.$  In this paper, we
consider ideal polytropic gas, that is, $P$ and $e$ satisfy \be  \la{1.3'}  P =R \theta/{v},\quad e=c_v\theta +\mbox{ const.},
\ee
where $R$ (specific gas constant) and $c_v $ (heat
capacity at constant volume) are both  positive constants. Moreover, for $\mu$ and $\ka,$
we  consider   the case where $\mu$ and $\ka$ are proportional to (possibly different) powers of  $\theta:$  \be   \la{1.3''}\mu=\tilde\mu \theta^\ga, \quad \kappa=\ti\kappa \theta^\beta, \ee
where $\ti\mu,\ti\ka>0$ and $\ga,\beta\ge 0$ are constants.

Then we impose the system \eqref{1.1}-\eqref{1.3''}  on
the initial condition
\be \la{1.4} (v(x,0),u(x,0),\theta(x,0))=(v_0(x),u_0(x),\theta_0(x)),\quad x\in
\O,\ee
and   three types of far-field and  boundary ones:

1) the Cauchy problem
\be \la{1.5}  \O=\r,\, \lim\limits_{|x|\rightarrow
	\infty}(v(x,t),u(x,t),\theta(x,t))=(1,0,1),\quad t>0;\ee

2)   boundary and far-field conditions for $ \O=(0,\infty),$
\be \la{1.6}  u(0,t)=0,\, \theta_x(0,t)=0,\,\lim\limits_{x\rightarrow
	\infty}(v(x,t),u(x,t),\theta(x,t))=(1,0,1),\quad t>0; \ee

3)   boundary and far-field conditions for $ \O=(0,\infty),$
\be \la{1.7} u(0,t)=0,\, \theta (0,t)=1,\, \lim\limits_{x\rightarrow
	\infty}(v(x,t),u(x,t),\theta(x,t))=(1,0,1), \quad t>0. \ee

According to the Chapman-Enskog expansion for the first level of approximation in kinetic theory, the viscosity $\mu$
and heat conductivity $\ka$ are functions of temperature alone (\cite{cip,cc}). Indeed, if the intermolecular
potential varies as $r^{-a},$  with intermolecular distance  $r$, then $\mu$ and $\ka$ are both
proportional to the power $(a+4)/(2a)$ of the temperature, that is, \eqref{1.3'} holds with $\ga=\beta=(a+4)/(2a).$  In  particular,  for elastic spheres $(a\rightarrow\infty),$ the dependence is like $\te^{1/2};$ while for Maxwellian molecules $(a=4),$ the dependence is linear.

For constant coefficients $(\ga=\b =0),$ Kazhikhov and Shelukhin
\cite{ks} first obtained    the global existence
of solutions  in bounded domains  for large initial data.  From then on,
significant progress has been made on the mathematical aspect of the initial and
initial
boundary value problems, see \cite{az1,az2,az3,akm,ji4,ji3,kaw} and the references therein. For the Cauchy problem  \eqref{1.1}-\eqref{1.5} and the initial boundary value
problems
\eqref{1.1}-\eqref{1.4}   \eqref{1.6}   and \eqref{1.1}-\eqref{1.4}   \eqref{1.7}
(in unbounded domains), Kazhikhov \cite{ka1} (also cf. \cite{akm,ji4}) obtains the global existence of strong solutions.

For the case that $\ga=0$ and $\beta\in(0,\infty),$  Pan-Zhang \cite{pz1} (see also \cite{jk1}) obtain the global strong solutions in bounded domains $\O$.  Recently, when $\O$ is unbounded,  for $\ga=0$ and $\beta>0,$ Li-Shu-Xu \cite{lsx1} obtained the global existence of the strong solutions which can be stated as follows:

\begin{lemma}[Global Existence{[23]}] \la{thm11}  Suppose that \be \ga=0,  \quad \beta> 0,\ee
	and  that the initial data $(v_0,u_0,\theta_0)$ satisfy
	\be\la{c5} v_0-1,u_0,\theta_0-1\in H^1(\O),\,\, \inf_{x\in \O}v_0(x)>0,\, \,
	\inf_{x\in \O}\theta_0(x)>0 ,\ee
	and   are compatible with \eqref{1.6},  \eqref{1.7}.
	Then there exists a unique global
	strong solution
	$(v ,u ,\theta )$   with positive  $v(x, t)$ and  $\theta(x, t)$ to the initial-boundary-value problem \eqref{1.1}-\eqref{1.5},  or \eqref{1.1}-\eqref{1.4},  \eqref{1.6}, or
	\eqref{1.1}-\eqref{1.4},  \eqref{1.7}   satisfying for any $T>0,$ \be
	\la{3k}\begin{cases}   v-1, \,u,\, \theta-1 \in L^\infty(0,T;H^1(\O) ), \\ v_t\in
		L^\infty(0,T;L^2 (\O))\cap L^2(0,T;H^1(\O)), \\ u_t,\,\theta_t,\,v_{xt},\,u_{xx},\,\theta_{xx} \,\in
		L^2(\O\times(0,T)).\end{cases}\ee
\end{lemma}

It is interesting to study the asymptotic behavior as  $t\rightarrow \infty $ of the solutions.   For $\ga=0$  and bounded domains, Kazhikhov \cite{kaz1} and Huang-Shi \cite{hs1} considered $ \beta=0$ and $  \beta>0$ respectively and obtained that  the global solution is asymptotically stable as time
tends
to infinity; see \cite{akm,az1,az2, na1,na2,na3,ni,qy} among others.
As far as the case of unbounded domains is concerned,
the asymptotic behavior as  $t\rightarrow \infty $ of the solution has been
studied under some smallness conditions on the initial data (see  \cite{ho,ji3,lz,ok,ka,knq,qy} and the references therein). For arbitrarily large initial data, only recently, Li-Liang \cite{ll} (see also \cite{ji1,ji2}) considered  the case that $\ga=0,\beta=0$ and established the large-time behavior of strong solutions. However,  it should
be mentioned here that the method  used in \cite{ll} relies heavily on the non-degeneracy of the heat
conductivity $\ka$ and cannot be applied directly to the degenerate and nonlinear case   $\beta>0.$ In this paper, we will prove that for $\beta>0$   the global strong solution  obtained by Lemma \ref{thm11} is asymptotically stable as time
tends to infinity   for large initial data. Our main result is as follows:

\begin{theorem} 
	
	\la{1k} Under the conditions of Lemma \ref{thm11}, let $(v ,u ,\theta
	)$ be the (unique) strong solution to \eqref{1.1}-\eqref{1.5}, or \eqref{1.1}-\eqref{1.4}, \eqref{1.6}, or \eqref{1.1}-\eqref{1.4},  \eqref{1.7} satisfying  \eqref{3k}  for any $T>0.$
	Then there exists a positive constant  $C_0$ depending only on
	$\ti\mu,\ti\ka,\beta, R,c_v,$ $\|(v_0-1,u_0,\theta_0-1)\|_{H^1(\O)},
	\inf\limits_{x\in \O}v_0(x),$ and $ \inf\limits_{x\in \O}\theta_0(x) $   such
	that
	\be \la{c2} C_0^{-1}\le  v (x,t)\le C_0,\,  C_0^{-1}\le \theta (x,t)\le C_0,\,\,\mbox{ for all } (x,t)\in \ol
	\O\times [0,\infty) ,\ee
	\be \la{c2a}\sup\limits_{0\le t<\infty} \|(v-1,u,\theta-1)\|_{H^1(\O)} +\int_0^\infty
	\left(\|v_x\|_{L^2(\O)}^2+\| (u_x,\theta_x)\|_{H^1(\O)}^2\right)dt\le C_0.\ee
	Moreover, the following large-time behavior holds
	\be\la{c3} \lim\limits_{t\rightarrow
		\infty}\left(\|(v-1,u,\theta-1)(t)\|_{L^p(\O)}+
	\|(v_x,u_x,\theta_x)(t)\|_{L^2(\O)}\right)=0,\ee for  any  $p\in (2,\infty].$
	
\end{theorem}

\begin{remark} Theorem \ref{thm11} can be regarded as a natural  generalization   of  Li-Liang's  result (\cite{ll}) where they considered the constant viscosity case ($\ga=\b=0$) to the degenerate and nonlinear one that $\ga=0,\b>0.$
	It is interesting to study the case that  $\ga>0,\b\ge 0,$  which will be left for future.
\end{remark}



We now comment on the analysis of this paper. Compared with  the constant viscosity case ($\ga=\b=0$) (\cite{ll}), the main difficulty comes from the   degeneracy and nonlinearity of the heat conductivity due to the fact that $\b>0.$ The key observations are as follows: First, stand calculations show that
for any $x\in \O  $ and $ N=[x],$ we have (see \eqref{DY})
\be\la{DYa}\ba v(x,t)  &= D_N(x,t) Y_N(t) \exp\left\{\int_0^t\frac{\te}{v}ds\right\}, \ea\ee
with \bnn\la{D}C^{-1}\leq D_N(x,t)\leq C,\enn and
\bnn \la{wq1.1} Y_N(t)\triangleq \exp\left\{\int_0^t\sigma (N,s)ds\right\}.\enn
Then, motivated by Huang-Shi \cite{hs1}, for $0\le \tau<t\le T,$    we have (see \eqref{ej1}) \bnn \la{ej1}\ba \int_\tau^t \sigma (N,s)ds  \le C-\frac{1}{2}\int_\tau^t \int_N^{N+1}\frac{\te}{v }dxds .
\ea\enn
This, after using the idea due to Jiang \cite{ji1,ji2}, implies
\bnn \la{wq1a2} Y_N(t)\le Ce^{-C^{-1}t},\quad \frac{Y_N(t)}{Y_N(\tau)}\le Ce^{-C^{-1}(t-\tau)},\enn which  combined  with \eqref{DYa} gives  the time-independently  lower and upper bounds of $v.$ Next, the second key step is to 
get the  time-independent $L^2$-norm   (in both space and time)  bound of $\theta^{-1/2}\theta_x$ (see   \eqref{w2.11} and \eqref{df8}). On the one hand, for $\beta\in(0,1),$ we just modified slightly the proof of \cite[Lemma 2.2]{ll}. Indeed,  since $v$ is bounded from above and below time-independently,   the standard   energetic  estimate (see \eqref{2.12})  shows that for $  \O_2(t) \triangleq\{x\in \O
|\,\theta(x,t)>2\} ,$  $$\int_0^\infty \int_{\O  \setminus\O_2(t)}\theta^\b \theta_x^2 dxdt  $$ is bounded. Hence, it suffices  to estimate the integral $$\mathcal{B}:= \int_0^\infty \int_{ \O_2(t)}\theta^\b  \theta_x^2dxdt .$$ In fact, motivated by
\cite{hxw,ll},  we multiply    the equation for the temperature (see \eqref{1.3}) by $(\theta-2)_+$ (see \eqref{1.1-1}). Then, to control the most difficult term  appearing in \eqref{1.1-1}, motivated by
\cite{hxw,ll}, we multiply  the equation for the velocity by $2u(\theta-2)_{+}$  (see \eqref{1.2-1}).    After some careful analysis on the integration  by parts over $\O_2(t)$ (see \eqref{df3}) and multiplying the equation for the velocity by $|u|^{2/(1-\b)} u,$ we finally find that $\mathcal{B}$ can be controlled by  (see \eqref{lia5}) $$\int_0^\infty \sup\limits_{x\in\O}(\theta-3/2)_+^2 (x,t)dt,$$
which in fact is bounded by $C(\ve)+C\ve \mathcal{B}  $ for any $\ve>0$ (see \eqref{lq22}). On the other hand, for $\b\in [1,\infty),$  multiplying the equation of the temperature (see \eqref{1.3}) by $ (\te^{-p}-4)_+$
with $p\ge 1,$ and after some careful calculations, we finally reach (see \eqref{ljj1})
\bnn   \int_0^T\int_{\O} \frac{\te^\b\te_x^2}{\te^{p+1}}dxdt+\int_0^T\int_{\O}\frac{u_x^2}{\te^p}dxdt  \le C(p),\enn which after choosing  $p=\b$ implies \bnn   \int_0^T\int_{\O} \te^{-1}\te_x^2  dxdt  \le C. \enn
Finally, to obtain   higher order estimates, we will modify some ideas due to \cite{lsx1} to obtain the estimates on the $L^2(\O\times(0,T))$-norms of both $u_t$ and $u_{xx}$ (see Lemmas  \ref{wl2.5} and \ref{lemm5}) which are crucial for further estimates on the upper bound of the temperature. These are the key to the proof of \eqref{c2}, and once \eqref{c2} is obtained, the proofs of \eqref{c2a} and \eqref{c3} are quite standard (see  \cite{ji3,lz,ok,ka, knq,qy}). The whole procedure will be carried out in the next section.

\section{Proof of Theorem \ref{1k}}

We begin with the following uniform (with respect to time) upper and lower bounds of $v.$

\begin{lemma} \la{lwq3} There exists a positive constant $C$ such that
	\be\ba\la{cqq1}
	C^{-1}\le v(x,t)\le C,
	\ea\ee
	where (and in what follows)   $C $ and  $C_i (i=1,2)$ denote
	generic positive constants
	depending only on $\ti\mu,\ti\ka,\beta,R,c_v,\|(v_0-1,u_0,\theta_0-1)\|_{H^1(\O)},
	\inf\limits_{x\in \O}v_0(x),$ and $ \inf\limits_{x\in \O}\theta_0(x).$\end{lemma}

\pf
First, using \eqref{1.1}, \eqref{1.2}, and \eqref{1.3'}, we rewrite
\eqref{1.'3} as \be  \la{1.3}  \theta_{t}+  \frac\theta v
u_{x}= \left(\frac{\te^\b\theta_{x}}{v}\right)_{x}+ \frac{u_{x}^{2}}{v}.
\ee  Multiplying  (\ref{1.1}) by
$ 1- {v}^{-1} $,  (\ref{1.2})  by $
u$,   (\ref{1.3}) by
$  1- {\theta}^{-1} $,  and adding them altogether, we obtain
\bess \ba&(u^2/2+ (v-\ln v-1)+ (\theta-\ln
\theta-1))_t+ \frac{u^2_x}{v\theta}
+ \frac{\te^\b\theta_x^2}{v\theta^2}\\&=\left(\frac{  u u_x}{v}-\frac{ u
	\theta}{v}\right)_x + u_x+ \left((1-\theta^{-1})
\frac{\theta^\b\theta_x}{v}\right)_x,\ea\eess which together with \eqref{1.5} or \eqref{1.6}
or \eqref{1.7} gives
\be\ba\la{2.12} &\sup_{0\le t\le T}\xix \left(  \frac{u^2}{2}+ (v-\log
v-1)+ (\theta-\log \theta-1)\right)dx +\int_0^TV(s)ds \le e_0,  \ea\ee
where \be V(s)\triangleq \xix \left(\frac{\theta^\b\theta_{x}^2}{v\te^2}+\frac{u_{x}^{2}}{v\te}
\right)(x,s)dx, \ee and
\bnn e_0\triangleq\xix\left(\frac{u_0^2}{2}+ (v_0-\log v_0-1)+ (\theta_0-\log
\theta_0-1)\right)dx.\enn

Then, for any $x\in \O  $ and $ N=[x],$  we have by \eqref{2.12}
\bnn\la{ja1}\int_{N}^{N+1}(v-\log v-1) + (\theta-\log\theta-1)dx\le e_0,\enn
which together with Jensen's inequality yields \be\la{ab} \al_1 \le \int_N^{N+1}v(x,t)dx\le \al_2,\quad\al_1\le \int_N^{N+1}\te(x,t)dx\le \al_2,\ee   where $0<\al_1<\al_2$  are two roots of  $$ y-\log y-1=e_0.$$ 

Next, letting $$\sigma\triangleq\frac{u_x}{v}-\frac{\te}{v}=(\log v)_t-\frac{\te}{v},$$ we   rewrite \eqref{1.2} as
\bnn\la{1.2'} u_t=\sigma_x.\enn
Integrating   this in $x$ over $[N,x] $ leads to
\be\la{w3.1} \int_N^xu_t(y,t)dy=\sigma(x,t)-\sigma(N,t),\ee which implies
\bnn  \int_N^x(u(y,t)- u_0(y))dy = \log v-\log v_0-\int_0^t\frac{\te}{v}ds-\int_0^t\sigma(N,s)ds.\enn This   gives
\be\la{DY}\ba v(x,t)  &= D_N(x,t) Y_N(t) \exp\left\{\int_0^t\frac{\te}{v}ds\right\}, \ea\ee
where
\bnn D_N(x,t)\triangleq v_0(x)\exp\left\{\int_{N}^x \left(u(y,t)-u_0(y)\right)dy\right\},
\enn and
\be \la{wq1.1} Y_N(t)\triangleq \exp\left\{\int_0^t\sigma (N,s)ds\right\}.\ee
Thus, it follows from  \eqref{DY} that
\be\la{DY'} v(x,t)=  D_N(x,t) Y_N(t) \left(1+ \int_0^t\frac{\te(x,\tau)}{ D_N(x,\tau) Y_N(\tau)}d\tau \right).\ee

Next, since \eqref{2.12} implies $$\left| \int_N^x  \left(u(y,t)-u_0(y)\right)dy\right|\leq \left(\dis\int_N^{N+1} u^2dy\right)^{1/2}+\left(\dis\int_N^{N+1} u_0^2dy\right)^{1/2}\leq C,$$  we have \be\la{D}C^{-1}\leq D_N(x,t)\leq C,\ee where and in what follows,  $C$ is a constant independent of $N$ and $T$.

Next, motivated by Huang-Shi \cite{hs1}, for $0\le \tau<t\le T,$ integrating \eqref{w3.1}   over $(N,N+1)\times (\tau,t)$ yields \be \la{ej1}\ba \int_\tau^t \sigma (N,s)ds &=\int_\tau^t\int_N^{N+1} \sigma (x,s)dx ds- \int_N^{N+1}\int_N^x u (y,t)dy dx\\& \quad+\int_N^{N+1}\int_N^x u (y,\tau)dy dx\\&\le C \int_\tau^t \int_N^{N+1}\frac{u_x^2}{v\te}dxds-\frac{1}{2}\int_\tau^t \int_N^{N+1}\frac{\te}{v }dxds+C\\&\le C-\frac{1}{2}\int_\tau^t \int_N^{N+1}\frac{\te}{v }dxds ,
\ea\ee where in the last inequality we have used \eqref{2.12}.

Then,  for \bnn\bar\te_N\triangleq\int_N^{N+1}\te(x,t)dx,\enn we obtain  from \eqref{ab} and Jensen's inequality that
\bnn\ba -\int_N^{N+1} \frac{\te}{v}dx&\le \int_N^{N+1} \left(-  \bar \te_N  +(\bar\te_N-\te)_+\right) \frac{1}{v}dx  \\  &\le  \left(-  \bar \te_N     + \max_{x\in[N,N+1]}\left( \bar\te_N -\te \right)_+\right)\left(\int_N^{N+1}vdx\right)^{-1} \\  &\le - 2C^{-1} + C\max_{x\in[N,N+1]}\left( \bar\te_N^{\b/2}-\te^{\b/2}\right)_+  \\&\le - 2 C^{-1}  + C\int_N^{N+1}\te^{-1+\b/2}|\te_x|dx \\&\le - C^{-1}  +C \int_N^{N+1}\frac{\te^{ \b } \te_x^2}{\te^2v}dx , \ea\enn
which together with \eqref{wq1.1},  \eqref{ej1},  and \eqref{2.12} gives
\be \la{wq1.2} Y_N(t)\le Ce^{-C^{-1}t},\quad \frac{Y_N(t)}{Y_N(\tau)}\le Ce^{-C^{-1}(t-\tau)}.\ee
Combining this, \eqref{DY'}, and \eqref{D}  gives
\be \la{wq1.3} v(x,t)\le Ce^{-C^{-1}t}+C\int_0^t \frac{Y_N(t)}{Y_N(\tau)}\te(x,\tau) d\tau. \ee

Next, it follows from \eqref{ab}  that
for all $(x,t)\in [N,N+1]\times [0,\infty),$
\bnn\ba\la{bqq2}  \left|\te^{\frac{\b+1}{2}}(x,t)- \bar\te_N^{\frac{\b+1}{2}}( t) \right| &\le \frac{\b+1}{2}\left(\int_N^{N+1} \frac{ \te^\b \te_x^2}{\te^2 v} dx\right)^{1/2}\left(\int_N^{N+1}  {\te v} dx\right)^{1/2}\\ &\le C V^{1/2}(t)  \max_{x\in [N,N+1]}v^{1/2}(x,t),\ea\enn which together with \eqref{ab} leads to
\be\ba\la{bqq7}
\frac{\alpha_1}{8}-CV(t)\max_{x\in [N,N+1]}v(x,t)\le  \theta(x,t)\le C +CV(t)\max_{x\in [N,N+1]}v(x,t).
\ea\ee
Putting this and \eqref{wq1.2} into \eqref{wq1.3} yields
\bnn\ba v(x,t)&\le C+C\int_0^t e^{-C^{-1}(t-\tau)}\left(1+V(\tau)\max_{x\in [N,N+1]}v(x,\tau)\right)d\tau\\ &\le C+C\int_0^t V(\tau)\max_{x\in [N,N+1]}v(x,\tau) d\tau,\ea\enn
which together with Gronwall's inequality and \eqref{2.12} shows that for any $t\in [0,\infty),$
\bnn \max_{x\in [N,N+1]}v(x,t)\le C.\enn
Since $x$ is arbitrary, we   have \be \la{wq2.20} v(x,t) \le C,\ee for all $(x,t)\in \O\times [0,\infty).$

Next, integrating \eqref{wq1.3} in $x$ over $(N,N+1),$ we get after using \eqref{ab} that\bnn \al_1\le Ce^{-C^{-1}t}+C\int_0^t\frac{Y_N(t)}{Y_N(\tau)} d\tau, \enn
which together with \eqref{DY'}, \eqref{D},   \eqref{wq1.2}, \eqref{bqq7}, and \eqref{wq2.20}   yields
\be  \la{wq2.23}\ba v(x,t) &\ge C^{-1}\int_0^t\frac{Y_N(t)}{Y_N(\tau)}\left(\frac{\al_1}{8} -CV(\tau)\right) d\tau \\&\ge C_1-Ce^{-C^{-1}t}-C\int_0^te^{-C^{-1}(t-\tau)} V(\tau) d\tau\\&\ge C_1 -Ce^{-C^{-1}t}-Ce^{-t/(2C)}\int_0^{t/2}  V(\tau) d\tau-C\int_{t/2}^t V(\tau) d\tau \\&\ge \frac{C_1}{2},\ea\ee provided $t\ge T_0 $ with $T_0$ independent of $T.$

Finally, it follows from \cite[Lemma 2.2]{lsx1} that there exists some positive constant $C$ independent of $T$ such that
\bnn v(x,t)\ge C^{-1}\enn for all $(x,t)\in \O\times [0,T_0].$  Combining this, \eqref{wq2.20}, and \eqref{wq2.23} gives \eqref{cqq1} and finishes the proof of Lemma \ref{lwq3}.
\thatsall


For further uses, we first state the following preliminary  estimates.
\begin{lemma} \la{lwq12} For any $p\ge 1,$   there exists some positive constant $C(p)$  such that
	\be \la{ljj1}  \int_0^T\int_{\O} \frac{\te^\b\te_x^2}{\te^{p+1}}dxdt+\int_0^T\int_{\O}\frac{u_x^2}{\te^p}dxdt  \le C(p).\ee \end{lemma}

\pf
First,   for $t\ge 0$ and $a>0,$ denote   \bess (\te>a) (t)\triangleq\{x\in \O
|\,\theta(x,t)>a\} ,\quad (\te<a) (t)\triangleq\{x\in \O
|\,\theta(x,t)<a\} .\eess For $\alpha>1,$ we derive from
\eqref{2.12} that\be \la{nep1}\sup\limits_{0\le t
	<\infty}\int_{(\te>\alpha)(t)}\theta dx\le C(\alpha)\sup\limits_{0\le t
	<\infty}\int_{\O}(\theta-\ln \theta-1)dx\le C(\alpha), \ee
and that  \be \la{nep2}\sup\limits_{0\le t <\infty}(|(\te>\alpha)(t)|+|(\te<\alpha^{-1})(t)|)\le C(\alpha).\ee

Next, it follows from \eqref{2.12} that \eqref{ljj1} holds for $p=1.$ Thus, it only remains to prove  \eqref{ljj1}   for    $p\ge  \b+4. $  Multiplying \eqref{1.3} by $ (\te^{-p}-4)_+\triangleq \max\{\te^{-p}-4,0\}$ and integrating the resulting equality over $\O\times (0,T)$ gives
\be\la{w2.6}\ba& \int_{\O} \int_{\te}^{4^{-1/p}}(s^{-p}-4)_+ dsdx +p\int_0^T\int_{(\te<4^{-1/p})(t)}v^{-1}\te^{\b-p-1}\te_x^2 dxdt \\&\quad+ \int_0^T\int_\O v^{-1} u_x^2(\te^{-p}-4)_+dxdt\\&=\int_{\O} \int_{\te_0}^{4^{-1/p}}(s^{-p}-4)_+ dsdx +\int_0^T\int_{\O}\te^{1-p}v^{-1}u_x
(1-4\te^p)_+ dxdt \\&\le C+ C(\ve)\int_0^T \int_{(\te<4^{-1/p})(t)} \te^{1-p} dx\max_{x\in \O}(1-4\te^p)_+^2dt\\&\quad + \ve \int_0^T\int_{\O}\te^{-p}u_x^2dxdt.\ea\ee

Next, direct computation yields that
\be \int_{\O} \int_{\te}^{4^{-1/p}}(s^{-p}-4)_+ dsdx\ge \frac{1}{p-1}\int_{(\te<4^{-1/p})(t)}\te^{1-p}dx-C ,\ee
and that \be\ba& \int_0^T\int_\O v^{-1} u_x^2(\te^{-p}-4)_+dxdt\\ &=\int_0^T\int_\O v^{-1} \te^{-p} u_x^2(1-4\te^{p})_+dxdt\\&\ge C^{-1}\int_0^T\int_{(\te<1/2)(t)}   \te^{-p} u_x^2dxdt\\&\ge C^{-1}\int_0^T\int_{\O}   \te^{-p} u_x^2dxdt-  C \int_0^T\int_{(\te>1/4)(t)}   \te^{ -1} u_x^2 dxdt\\&\ge C^{-1}\int_0^T\int_{\O}   \te^{-p} u_x^2dxdt-  C.\ea\ee

Moreover, it follows from \eqref{2.12} and \eqref{nep2} that
\be\ba\int_0^T\max_{x\in \O}(1-4\te^p)_+^2dt &=\int_0^T\left(\max_{x\in \O}\int_x^{ \infty} (-(1-4\te^p)_+)_ydy\right)^2dt \\&\le C\int_0^T\left(\int_{(\te<4^{-1/p})(t)}
\te^{-1+\b/2}|\te_y|dy\right)^2dt\\&\le C\int_0^T \int_{\O}
\te^{ \b-2} \te_y^2dydt\le C, \ea\ee and that \be\la{w2.10}\ba &\int_0^T\int_{\O} \te^{\b-p-1}\te_x^2 dxdt \\&\le C \int_0^T\int_{(\te<4^{-1/p})(t)}v^{-1}\te^{\b-p-1}\te_x^2 dxdt+C\int_0^T\int_{(\te>1/2 )(t)} \te^{\b-2}\te_x^2 dxdt \\&\le C \int_0^T\int_{(\te<4^{-1/p})(t)}v^{-1}\te^{\b-p-1}\te_x^2 dxdt+C. \ea\ee
Combining \eqref{w2.6}--\eqref{w2.10} and Gronwall's inequality gives \eqref{ljj1}  and finishes the proof of Lemma \ref{lwq12}. \thatsall

Next, the following key Lemmas  \ref{wl2.1} and \ref{l3.2}  will deal with the  $L^2$-norm (in both space and time) bound  of both $ u_x$ and $\theta^{-1/2}\theta_x.$

\begin{lemma} \label{l3.2} For $\b\in(0,1),$ there exists some positive constant $C $   such that for
	any $T>0,$
	\be  \la{df8}
	\int_{0}^{T}\int_{\O}\left( u_{x}^{2}+
	(\te^{-1}+\te^\b)\theta_{x}^{2}\right)dxdt
	\leq C. \ee
\end{lemma}

\pf The proof of Lemma \ref{l3.2} will be divided into three steps.

{\it Step 1.}
First,    integrating  \eqref{1.3} multiplied  by $(\theta-2)_{+}\triangleq \max\{\theta-2,0\}$ over $\O\times(0,T)$
gives
\begin{equation}\begin{split} \label{1.1-1}  &\frac{ 1}{2} \int_{\O}
		(\theta-2)_{+}^{2}dx-\frac{1 }{2}\int_{\O} (\theta_{0}-2)_{+}^{2}dx+ \int_{0}^{T}\int_{(\te>2)(t)}\frac{ \te^\b \theta_x^{2}}{v}dxdt
		\\  & =-   \int_{0}^{T}\int_{\O}\frac{\theta}{v}u_{x}(\theta-2)_{+}dxdt + \int_{0}^{T}\int_{\O}
		\frac{u_{x}^{2}}{v}(\theta-2)_{+}dxdt. \end{split}\end{equation}

To estimate the last term on the right hand side of \eqref{1.1-1}, motivated by
\cite{ll},  we multiply   \eqref{1.2} by $2u(\theta-2)_{+}$ and integrate the
resulting equality over $\O\times (0,T)$ to get
\be \ba\begin{split} \label{1.2-1}&  \int_{\O}  u^2  (\theta-2)_{+}dx
	-\int_{\O}   u_0^2
	(\theta_0-2)_{+}dx+2 \int_{0}^{T}\int_{\O} \frac{u_{x}^{2}}{v}(\theta-2)_{+}dxdt\\&= 2  \int_{0}^{T}\int_{\O} \frac{\theta}{v} u_{x}(\theta-2)_{+}dxdt
	+ 2  \int_{0}^{T}\int_{(\te>2)(t)}\frac{\theta}{v} u \theta_xdxdt\\&\quad
	- 2 \int_{0}^{T}\int_{(\te>2)(t)}\frac{u_{x}}{v}u  \theta_xdxdt
	+\int_{0}^{T}\int_{(\te>2)(t)} u^{2} \theta_tdxdt. \end{split}\ea\ee
Adding  (\ref{1.2-1}) to (\ref{1.1-1}), we obtain after using \eqref{1.3} that
\be \begin{split} \label{1.3-1}
	&\int_{\O}\left[\frac{1 }{2} (\theta-2)_{+}^{2}+  u^{2}(\theta-2)_{+}\right]dx-\int_{\O}\left[\frac{1 }{2} (\theta_0-2)_{+}^{2}+
	u_0^{2}(\theta_0-2)_{+}\right]dx\\&\quad+
	\int_{0}^{T}\int_{(\te>2)(t)}\frac{\te^\b\theta_x^{2}}{v}dxdt
	+ \int_{0}^{T}\int_{\O}\frac{u_{x}^{2}}{v}(\theta-2)_{+}dxdt
	\\  &=  \int_{0}^{T}\int_{\O}
	\frac{\theta}{v}u_{x}(\theta-2)_{+}dxdt
	+ 2 \int_{0}^{T}\int_{(\te>2)(t)} \frac{\theta}{v} u \theta_xdxdt\\&\quad
	- 2 \int_{0}^{T}\int_{(\te>2)(t)}\frac{u_{x}}{v}u  \theta_x dxdt +  \int_{0}^{T}\int_{(\te>2)(t)}u^{2}
	\left( \frac{u_{x}^{2}}{v}- \frac{\theta}{v}u_{x}\right)dxdt\\
	&\quad +  \int_{0}^{T}\int_{(\te>2)(t)}u^{2}
	\left(\frac{\te^\b\theta_{x}}{v}\right)_{x}dxdt \triangleq \sum_{i=1}^{5}I_{i}. \end{split}\ee

We estimate each $I_i(i=1,\cdots,5)$ as follows:

First, it follows from Cauchy's inequality
that\be \ba\la{df1} |I_1|&=  \left|\int_{0}^{T}\int_{\O}
\frac{\theta}{v}u_{x}(\theta-2)_{+}dxdt \right|\\
&\le
\frac{1}{2}\int_{0}^{T}\int_{\O}\frac{u_{x}^{2}}{v}(\theta-2)_{+}dxdt +C\int_{0}^{T} \int_{\O}\theta^{2}(\theta-2)_{+}dxdt\\
&\le  \frac{1}{2}\int_{0}^{T}\int_{\O}\frac{u_{x}^{2}}{v}(\theta-2)_{+} dxdt +C\int_{0}^{T}\int_{\O } \theta(\theta-3/2)^{2}_{+}dxdt\\
&\le  \frac{1}{2}\int_{0}^{T}\int_{\O}\frac{u_{x}^{2}}{v}(\theta-2)_{+}dxdt +C\int_{0}^{T}\sup_{x\in
	\Omega }(\theta-3/2)^{2}_{+}(x,t)dt,\ea\ee   where in the last inequality we have used  \eqref{nep1}.

Next, Cauchy's inequality     yields that for any $\ve>0,$  \be \ba\la{df2}
|I_{2}|+|I_3|&=2 \left|\int_{0}^{T}\int_{(\te>2)(t)}
\frac{\theta}{v}u\theta_xdxdt\right|+2
\left|\int_{0}^{T}\int_{(\te>2)(t)}\frac{u_{x}}{v}u \theta_xdxdt\right|\\
&\le C\int_{0}^{T}\int_{\O}\theta_x^2dxdt+C \int_{0}^{T}
\int_{(\te>2)(t)}u^{2}\theta^{2}dxdt +C \int_{0}^{T}\int_{\O}u^{2}u_{x}^{2} dxdt \\
&\le C(\ve)+\ve\int_{0}^{T}\int_{\O}\te^\b \theta_x^2dxdt
+C \int_{0}^{T}\sup_{x\in
	\Omega }(\theta-3/2)^{2}_{+}(x,t)dt\\&\quad +C \int_{0}^{T}
\int_{\O} u_{x}^{2} dxdt+C \int_{0}^{T}
\int_{\O}|u|^{2/(1-\b)} u_{x}^{2}dxdt , \ea\ee where in the last inequality we have used\be  \la{2.13}\ba\int_{0}^{T}
\int_{(\te>2)(t)}u^{2}\theta^{2}dxdt &\le 16\int_{0}^{T}
\int_{\Omega }u^{2}(\theta-3/2)^{2}_{+}dxdt \\&\le C  \int_{0}^{T}\sup_{x\in
	\Omega }(\theta-3/2)^{2}_{+}(x,t)dt,\ea \ee due to \eqref{2.12}.

Then,  it follows from Cauchy's inequality  and \eqref{2.13} that
\be\la{df}\ba
|I_4| &\le  C \int_{0}^{T}
\int_{\O} u_{x}^{2}dxdt +C \int_{0}^{T}
\int_{\O}|u|^{2/(1-\b)} u_{x}^{2}dxdt\\&\quad+ C \int_{0}^{T}\sup_{x\in
	\Omega }(\theta-3/2)^{2}_{+}(x,t) dt.
\ea\ee

Finally, for $\eta>0  $ and
\bnn \ba   \varphi_{\eta}(\theta)\triangleq
\begin{cases}
	1, & \theta-2>\eta, \\
	(\theta-2)/\eta, & 0\le \theta-2\le \eta,\\
	0, &  \theta-2\le 0,
\end{cases}
\ea\enn
Lebesgue's dominated convergence theorem shows that for $\b<1$ and any $\ve>0,$ \be \ba\la{df3} I_{5} &=
\lim_{\eta\rightarrow0+}\int_{0}^{T} \int_{\O}
\varphi_{\eta} (\theta) u^{2}\left(\frac{\te^\b\theta_{x}}{v}\right)_{x} dxdt \\
&= \lim_{\eta\rightarrow0+} \int_{0}^{T}\int_{\O} \left(-2 \varphi_{\eta}(\theta ) uu_{x}\frac{\te^\b\theta_{x}}{v}-
\varphi_{\eta}' (\theta ) u^{2}\frac{\te^\b\theta_{x}^{2}}{v} \right)dxdt \\
&\le - {2 }   \int_{0}^{T}\int_{(\te>2)(t)}   uu_{x}\frac{\te^\b\theta_{x}}{v} dxdt\\
&\le \ve\int_{0}^{T}\int_{\O}  \left( u_{x}^{2}\te +\te^\b\theta_x^{2} \right)dxdt
+C(\ve)\int_{0}^{T}\int_{\O}  |u|^{2/(1-\b)}u_{x}^{2}dxdt ,
\ea\ee where in the third inequality we have used
$\varphi_\eta'(\theta)\ge 0.$

Noticing that
\bess\ba  &\int_{0}^{T}\int_{\O}\left(u_{x}^{2}\theta+\te^\b\theta_{x}^{2}\right)dxdt\\
&\le\int_{0}^{T}\int_{(\te>3)(t)}\left(u_{x}^{2}\theta +\te^\b\theta_{x}^{2}\right)dxdt+ \int_{0}^{T}\int_{(\te<4)(t)}\left(u_{x}^{2}\theta +\te^\b\theta_{x}^{2}\right)dxdt \\
&\le C\int_{0}^{T}\int_{(\te>2)(t)}\frac{1}{v}\left(  u_{x}^{2}(\theta-2)_{+}+ \te^\b
\theta_x^{2}\right)dxdt
+C,\ea\eess where in the second inequality we have used \eqref{cqq1}  and \eqref{2.12},
we substitute \eqref{df1}, \eqref{df2}, \eqref{df}, and \eqref{df3} into \eqref{1.3-1} and choose $\ve$ suitably
small to obtain
\be \begin{split} \label{df4}
	&\sup_{0\le t\le T}\int_{\O}
	(\theta-2)_{+}^{2}dx+\int_{0}^{T}\int_{\O} \left(u_{x}^{2}\theta+\te^\b\theta_{x}^{2}\right)dxdt\\
	&\le  C+ C  \int_{0}^{T}\sup_{x\in
		\Omega }(\theta-3/2)^{2}_{+}(x,t)dt +
	C_2\int_{0}^{T}\int_{\O}|u|^{2/(1-\b)}u_{x}^{2}dxdt, \end{split}\ee
where   we have used the following simple fact that for any $\de>0,$\be\la{1pp} \ba 2\int_{0}^{T}\int_{\O}u_x^2dxdt&\le \de\int_{0}^{T}\int_{\O}\theta u_x^2 dxdt +  \de^{-1} \int_{0}^{T}\int_{\O}\theta^{-1} u_x^2 dxdt\\&\le  \de \int_{0}^{T}\int_{\O}\theta u_x^2 dxdt +C (\de),\ea \ee due to Cauchy's inequality, \eqref{2.12}, and \eqref{cqq1}.

{\it Step 2.} To estimate the last term on the right hand side of \eqref{df4}, we multiply
\eqref{1.2} by $|u|^{\a }u (\a= 2/(1-\b))$ and integrate the resulting equality over $\O\times (0,T)$
to get
\be \ba\la{lia}  &\frac{1}{\a+2}\int_{\O}  |u|^{\a+2}dx+(\a+1)
\int_{0}^{T}\int_{\O} \frac{|u|^{\a}u_{x}^{2}}{v} dxdt\\
&\le C+C\int_{0}^{T}\int_{\O} \left(\frac{ |1-v|
}{v}+ \frac{|\theta-1|1_{(\te<3)(t)}}{v} \right)|u|^{\a }|u_{x}|dxdt \\&\quad+C\int_{0}^{T}\int_{(\te>2)(t)} \frac{|\theta-1|}{v} |u|^{\a}|u_{x}|dxdt \triangleq C+\sum\limits_{i=1}^2 J_i.\ea\ee

It follows from \eqref{2.12}
and \eqref{cqq1} that for any $\de\in [2,4],$\be
\la{4y} \ba&\sup\limits_{0\le t\le T}\int_\O  (v-1)^2dx+\sup_{0\le t\le T}\int_{(\te<\de)(t)}(\theta-1)^{2}dx \\&\le  C\sup_{0\le t\le T}\int_{\O}(v-\ln v-1)dx+ C \sup_{0\le t\le T}\int_{\O}(\theta-\ln \theta-1)dx \le C ,\ea\ee
which together with Holder's inequality yields that
\be\la{cc}\ba |J_1| &\le C\int_0^T  \||u|^\a u_x \|_{L^2(\O)}\left(\int_\O  (v-1)^2dx+ \int_{(\te<3)(t)}(\theta-1)^{2}dx\right)^{1/2}\\&\le C\int_0^T \||u|^\a u_x \|_{L^2(\O)} . \ea\ee

Next,  on the one hand, if $\b\le 1/2,$
\be \la{w2.31}\ba \||u|^\a u_x \|_{L^2(\O)} &\le C \max\limits_{x\in\O}  u^2(x,t) \||u|^{\a\b}u_x \|_{L^2(\O)} \\&\le C(\ve)\int_\O u_x^2dx +\ve \int_\O |u|^{\a } u_x^2dx , \ea\ee  where in the second inequality we have used   \eqref{2.12} and the following simple fact that for any $w\in H^1(\O),$
\be  \la{pj}\ba \max\limits_{x\in\O} w^2(x ) &=\max\limits_{x\in\O}\left(-2 \int_x^\infty w(y )w_x(y )dy\right)\\&\le 2\|w \|_{L^2(\O)}\|w_x \|_{L^2(\O)}. \ea\ee On the other hand, if $\b\in (1/2,1),$ we have
\be \la{w2.33}\ba \||u|^\a u_x \|_{L^2(\O)} &\le C(\ve) \max\limits_{x\in\O}  |u|^\a(x,t)+\ve \int_\O |u|^{\a } u_x^2 dx\\&\le C(\ve)\int_\O u_x^2 dx+2\ve \int_\O |u|^{\a } u_x^2 dx, \ea\ee
where we have used
\bnn  \ba
\max\limits_{x\in\O}|u|^{\a }
&= \max\limits_{x\in\O}\int^{\infty}_{x}(-  \pa_x |u|^{\a } )dx\\
&\le C    \int_\O  |u_x| |u|^{\a-1}dx\\
&\le C   \left( \int_\O  |u|^{\a-4 } u_x^2 dx \right)^{1/2} \left(\int_\O   |u|^{\a+2 }dx\right)^{1/2}\\
&\le C\left( \int_\O  |u|^{\a -4} u_x^2 dx \right)^{1/2}
\max\limits_{x\in\O}|u|^{  \a/2}, \ea\enn
due to  $\a=2/(1-\b)>4.$ Thus,  combining \eqref{cc}, \eqref{w2.31}, and \eqref{w2.33} implies for   $\b\in(0,1),$
\be \la{cc'}\ba |J_1|\le   C(\ve)\int_0^T  \int_\O u_x^2dx dt+C\ve \int_0^T\int_\O |u|^{\a } u_x^2dxdt . \ea\ee

Next, it follows from   Cauchy's inequality and \eqref{2.13} that
\be \ba\la{rf1}
|J_2|&\le   C\int_{0}^{T}\int_{(\te>2)(t)}   \theta^{2}|u|^{\a}dxdt+   \int_{0}^{T}\int_{\O}   \frac{|u|^{\a}u_{x}^{2}}{v}dxdt\\
&\le  C\int_{0}^{T}
\max_{x\in
	\Omega }(\theta-7/4)^{\frac{2\a}{\a+2}}_{+} \int((\te-1)^2+|u|^{\a+2} )dxdt\\&\quad+  \int_{0}^{T}\int_{\O}   \frac{|u|^{\a}u_{x}^{2}}{v} dxdt \\
&\le C\int_{0}^{T}
\max\limits_{x\in\O}((\theta(x,t)-3/2)_{+}^{\b+2}\te^{-1}) \int((\te-1)^2+|u|^{\a+2} )dxdt\\&\quad+   \int_{0}^{T}\int_{\O}   \frac{|u|^{\a}u_{x}^{2}}{v}dxdt ,
\ea\ee
due to $\b\in (0,1) $ and $2\a/(\a+2)<\b+1.$

Direct computation gives
\bnn \ba 
& \max\limits_{x\in\O}((\theta(x,t)-3/2)_{+}^{\b+2}\te^{-1})\\
&= \max\limits_{x\in\O}\int^{\infty}_{x}(-  \pa_x( (\theta-3/2)_{+}^{\b+2} \te^{-1} ))dx\\
&\le C    \int_\O  |\theta_x| (\theta-3/2)_{+}^{\b+1} \te^{-1}dx+ C    \int_\O  |\theta_x| (\theta-3/2)_{+}^{\b+2} \te^{-2}dx  \\
&\le C   \left( \int_\O  \theta_x^2  \theta^{\b-2 }dx\right)^{1/2} \max\limits_{x\in\O}((\theta(x,t)-3/2)_{+}^{\b+2}\te^{-1})^{1/2}   ,  \ea\enn where in the last inequality we have used \eqref{nep1}. This implies \be  \la{w2.37} \ba
\max\limits_{x\in\O}((\theta(x,t)-3/2)_{+}^{\b+2}\te^{-1}) \le C \int_\O  \theta_x^2 \theta^{\b-2 }dx. \ea\ee

Then, putting   \eqref{cc'},  \eqref{rf1},  and \eqref{w2.37} into \eqref{lia} and choosing $\ve$ suitably small  gives
\be \ba\la{1.5-1}
& \sup_{0\le t\le T}\int_{\O} | u|^{\a+2}dx+\int_{0}^{T}\int_{\O} |u|^{\a}u_{x}^{2}dxdt\\
&\le C+C \int_{0}^{T}\int_{\O}
\te^{\b-2}\te_x^2dx \int_{\O}((\te-1)^2+|u|^{\a+2} )dxdt+C  \int_{0}^{T}\int_{\O} u_{x}^{2}dxdt \\
&\le C \int_{0}^{T}\int_{\O}
\te^{\b-2}\te_x^2dx \int_{\O}((\te-1)^2+|u|^{\a+2} )dxdt + C(\de)\\&\quad +C\delta \int_{0}^{T}\int_{\O} \theta u_{x}^{2} dxdt,
\ea\ee where in the last inequality we have used \eqref{1pp}.

Adding \eqref{1.5-1}  multiplied  by $  C_2+1$ to   \eqref{df4}, then choosing
$\delta $ suitably  small, we have
\be \ba\label{lia5}
&\sup_{0\le t\le T}\int_{\O} \left[(\theta-1)^{2}+ |u|^{\a+2}\right]dx+
\int_{0}^{T}\int_{\O}\left[(\theta+|u|^{\a})u_{x}^{2}+
\te^\b\theta_{x}^{2}\right]dxdt\\
&\leq  C  +C \int_{0}^{T}\int_{\O}
\te^{\b-2}\te_x^2 dx\int_{\O}((\te-1)^2+|u|^{\a+2} )dxdt\\&\quad +C \int_{0}^{T}  \sup_{x\in
	\Omega }(\theta-3/2)^{2}_{+}(x,t)dt ,
\ea\ee where we have used \bnn \ba\int_\O (\te-1)^2dx&\le C \int_{(\te>3)(t) } (\te-1)^2dx+ \int_{(\te<4)(t) }(\te-1)^2dx \\&\le C \int_{\O } (\te-2)_+^2dx+C,\ea\enn due to \eqref{4y}.

{\it Step 3.} It remains to estimate the last   term  on the right hand side of
\eqref{lia5}.
Indeed, choosing $\de=-1$ in \eqref{lia22}  yields   that for    any $\ve>0,$
\be \la{lq22} \ba
\int_0^T\sup\limits_{x\in\O}(\theta(x,t)-3/2)_{+}^{2}dt  & \le C(\ve)\int_0^T
\int_{\O}\frac{\te^\b\theta_{x}^{2}}{v\theta^2}dxdt+ \ve \int_0^T \int_{\O} {\te^\b\theta_{x}^{2}}dxdt\\ & \le C(\ve) + \ve \int_0^T \int_{\O} {\te^\b\theta_{x}^{2}}dxdt, \ea\ee due to $\b<1,$ \eqref{nep1}, and  \eqref{2.12}.
Putting this   into  \eqref{lia5}, choosing $\ve$ suitably small, and using   Gronwall's inequality  lead  to
\bess  \sup_{0\le t\le T}\int_{\O} \left[(\theta-1)^{2}+ |u|^{\a+2}\right]dx+
\int_{0}^{T}\int_{\O}\left[(\theta+|u|^{\a})u_{x}^{2}+
\te^\b\theta_{x}^{2}\right] dxdt\leq   C ,
\eess   which combined with   \eqref{1pp}   and \eqref{2.12}  immediately gives \eqref{df8}. The proof of Lemma \ref{l3.2} is
completed. \thatsall

\begin{lemma} \la{wl2.1} For $\b\in[1,\infty),$ there exists some positive constant $C $   such that for
	any $T>0,$
	\be\la{w2.11}
	\int_{0}^{T}\int_{\O} ( u_{x}^{2}+\te^{-1}\te_x^2 )dxdt
	\leq C. \ee
\end{lemma}

\pf 
First,  choosing $p=\b$  in \eqref{ljj1}  gives \bnn \la{ljj1a} \int_{0}^{T}\int_{\O} \te^{-1}\te_x^2 dxdt
\leq C.\enn

Then, it remains to prove \be \la{ljj2a} \int_{0}^{T}\int_{\O}  u_x^2 dxdt
\leq C.\ee   Indeed,     multiplying
\eqref{1.3} by $  (\te-2)_+\te^{-1}$ and integrating the resulting equality over $\O\times(0,T)$ yields
\be \la{w2.13}\ba  & \int_{0}^{T}\int_{\O}
\frac{u_{x}^{2}}{v}(\theta-2)_{+}\te^{-1}  dxdt\\  & =2 \int_{0}^{T}\int_{(\te>2)(t)}\frac{ \te^{\b-2} \theta_x^{2}}{v}dxdt
+  \int_{\O  }\int_2^\te (s-2)_+s^{-1}dsdx\\&\quad -  \int_{\O  }\int_2^{\te_0} (s-2)_+s^{-1}ds dx +    \int_{0}^{T}\int_{\O}\frac{ (\theta-2)_{+}}{v}u_{x}dxdt\\&\le C+C(\ve) \int_0^T\sup_{x\in
	\Omega }(\theta-3/2)^{\b+1}_{+}(x,t)dt
+\ve \int_0^T\int_\O u_x^2 dxdt . \ea\ee

Then, direct calculation shows that for $\de\ge -1,$
\bnn \ba    &
\sup\limits_{x\in\O}(\theta(x,t)-3/2)_{+}^{\de+3}\\
&= \sup\limits_{x\in\O}\int^{\infty}_{x}(-  \pa_x(\theta-3/2)_{+}^{\de+3}  )dx\\
&\le C    \int_\O  |\theta_x| (\theta-3/2)_{+}^{\de+2}dx \\
&\le C   \left( \int_{(\te> 3/2)(t)}  \theta_x^2  \theta^{\de }dx\right)^{1/2} \left(\int_{(\te> 3/2)(t)} (\theta(x,t)-3/2)_{+}^{2\de+4}   \theta^{-\de} dx\right)^{1/2} \\
&\le C   \left( \int_{(\te> 3/2)(t)}  \theta_x^2  \theta^{\de }dx\right)^{1/2} \sup\limits_{x\in\O}(\theta(x,t)-3/2)_{+}^{(\de+3) / 2} ,\ea\enn where in the last inequality we have used \eqref{nep1}. This
gives
\be    \la{lia22} \ba
\sup\limits_{x\in\O}(\theta(x,t)-3/2)_{+}^{\de+3}  \le C   \int_{(\te> 3/2)(t)}  \theta_x^2 \theta^{\de } dx . \ea\ee
In particular, since $\b\ge 1,$ choosing $\de=\b -2$ in \eqref{lia22} gives \be  \la{w2.15} \int_0^T\sup\limits_{x\in\O}(\theta(x,t)-3/2)_{+}^{\b+1} dt\le C.\ee

Finally, it follows from \eqref{2.12}, \eqref{w2.13}, and \eqref{w2.15} that \bnn\ba \int_0^T\int_\O u_x^2dxdt&\le C\int_0^T\int_{(\te>3)(t)} \frac{u_{x}^{2}}{v}(\theta-2)_{+}\te^{-1}dxdt\\ &\quad+C \int_0^T\int_{(\te<4)(t)} \frac{u_{x}^{2}}{v\te}dxdt\\&  \le C(\ve)+C\ve \int_0^T\int_\O u_x^2dxdt,\ea\enn
which   gives \eqref{ljj2a} and finishes the proof of Lemma \ref{wl2.1}. \thatsall

Next, we will  derive some necessary  uniform estimates on the spatial derivatives  of   $v.$
\begin{lemma}\la{wl2.5}    There exists some positive constant $C $   such that for any
	$T>0,$ \be\la{z4} \sup_{0\le t\le T}\int_{\O} v_x^2dx
	+\int_0^T\int_{\O}\left(1+\theta \right) v_x^2dxdt\le
	C.\ee\end{lemma}

{\it Proof.} First, integrating  \eqref{1.2}     multiplied by  $\frac{v_{x}}{v}$
over $\O ,$ we obtain after  using \eqref{1.1}  that
\be\la{w2.43} \begin{split} &\frac{1}{2}\frac{d}{dt}\int_{\O} \frac{v_{x}^2}{v^2}dx+\frac{1}{2} \int_{\O}
	\frac{  v^{2}_{x}}{v^{3}}dx +\frac{1}{2} \int_{\O}
	\frac{\theta v^{2}_{x}}{v^{3}}dx
	\\ &=   \int_{\O} \left(\frac{\theta}{v}\right)_{x}\frac{v_{x}}{v}dx+ \int_{\O} u_t
	\frac{v_{x}}{v}dx +\frac{1}{2} \int_{\O}
	\frac{  v^{2}_{x}}{v^{3}}dx +\frac{1}{2} \int_{\O}
	\frac{\theta v^{2}_{x}}{v^{3}}dx \\&=  \int_{\O} \frac{\theta_{x}v_{x}}{v^{2}}dx+\frac{1}{2} \int_{\O}
	\frac{(1-\theta) v^{2}_{x}}{v^{3}}dx + \frac{d}{dt}\int_{\O} u \frac{v_{x}}{v}dx+\int_{\O}
	u_x\frac{ v_t}{v}dx\\&\le   C \int_{\O} \frac{\theta_{x}^2 }{\te}dx+\frac{1}{8} \int_{\O}
	\frac{\theta v^{2}_{x}}{v^{3}}dx+\frac{1}{8} \int_{\O}
	\frac{ v^{2}_{x}}{v^{3}}dx+C \max_{x\in\O}(1-\theta)_+^4\int_{\O}
	\frac{ v^{2}_{x}}{v^{2}}dx \\&\quad+ \frac{d}{dt}\int_{\O} u \frac{v_{x}}{v}dx+\int_{\O}
	\frac{u_x^2}{v}dx.\end{split}\ee

Then, it follows from \eqref{w2.11}, \eqref{df8}, and \eqref{4y} that
\be\la{w2.44}\ba   \int_0^T\max_{x\in\O}(1-
\theta)^4_+dt   &= 4\int_0^T\max_{x\in\O} \left(\int_x^{\infty}(1-
\theta)_+ \te_xdx\right)^2dt \\&\le C\int_0^T\left(\int_{\O} (1-
\theta)_+ \te^{-1/2}|\te_x|dx\right)^2dt  \\&\le C \int_0^T\int_{\O} (1-
\theta)_+^2 dx  \int_{\O} \te^{-1}\te_x^2dxdt\le C. \ea\ee
Applying Gronwall's inequality to  \eqref{w2.43},  we obtain \eqref{z4} after using  \eqref{w2.44}, \eqref{w2.11}, and \eqref{df8}.  The proof of Lemma \ref{wl2.5} is
completed. \thatsall

\begin{lemma}\la{lemm5} There exists a   positive constant $C $ such that  for any $T\ge 0,$
	\be\ba\la{lm5}  \sup_{0\le t\le T}\int_{\O} u_x^2dx+\int_0^T\int_{\O}(u_t^2+u_x^2+ u_{xx}^2+\te_x^2)dxdt\le C.\ea\ee\end{lemma}

\pf First, we rewrite the momentum equation  \eqref{1.2}  as
\be\ba\la{mom5} u_t-\frac{u_{xx}}{v}=-\frac{u_xv_x}{v^2}-\frac{\te_x}{v } +\frac{\te v_x}{v^2}.\ea\ee
Multiplying both sides of \eqref{mom5} by $u_{xx}$ and integrating the resultant equality in $x$ over $\O$ lead to
\be\la{qu2}\ba
&\frac{1}{2}\frac{d}{dt}\int_{\O} u_x^2dx+\int_{\O}\frac{u_{xx}^2}{v} dx\\
&\le \left|\int_{\O}\frac{u_xv_x}{v^2}u_{xx}dx\right|+\left|\int_{\O}\frac{\te_x}{v }u_{xx}dx\right| +\left|\int_{\O}\frac{\te v_x}{v^2}u_{xx}dx\right|\\
&\le \frac{1}{4}\int_{\O}\frac{u_{xx}^2}{v} dx+C\int_{\O}\left(u_x^2v_x^2+v_x^2\te^2 +\te_x^2\right)dx.
\ea\ee

Then, direct computation yields that for any $\de>0,$
\be\la{qu1}\ba &\int_{\O} \left(u_x^2v_x^2+v_x^2\te^2+\te_x^2\right)dx\\ &\le C\left(\max_{x\in\O} u_x^2+\max_{x\in[0,1]} (\te-3/2 )_+^2+1\right)\int_{\O}  v_x^2dx +\int_{\O} \te_x^2dx \\ &\le \eta \int_{\O} u_{xx}^2dx+ C(\eta) \int_{\O} ( u_x^2+v_x^2+\te_x^2)dx+C\max_{x\in\O}\left(\te-3/2\right)_+^2  ,\ea\ee
where in the last inequality  we have used
\be\ba\la{ux2}
\max_{x\in\O}u_x^2 \le \eta \int_{\O} u_{xx}^2dx+ C(\eta) \int_{\O} u_x^2dx,
\ea\ee due to \eqref{pj}. Thus, noticing that combining \eqref{lia22} where we choose $\de=-1$, \eqref{w2.11}, and \eqref{df8} gives \be\la{vx1} \int_0^T\max_{x\in\O} (\te-3/2 )_+^2dt\le C,\ee
putting \eqref{qu1} into \eqref{qu2} and choosing $\eta$ suitably small, we have
\be\ba\la{le5eq1}\sup_{0\le t\le T}  \int_{\O} u_x^2dx +\int_0^T\int_{\O} u_{xx}^2 dxdt\le C   +C\int_0^T \int_{\O}\te_x^2dx dt,\ea\ee due to \eqref{w2.11}, \eqref{df8}, and \eqref{z4}.

Next, on the one hand, if $\b\ge 2,$ choosing $p=\b-1$ in \eqref{ljj1} gives
\bnn\ba\la{sq1} \int_0^T\int_{\O}  \te_x^2dxdt \le C, \ea\enn
which  along with    \eqref{le5eq1}  shows
\be\ba\la{lm5w}  \sup_{0\le t\le T}\int_{\O} u_x^2dx+\int_0^T\int_{\O} u_{xx}^2dxdt+\int_0^T\int_{\O}\te_x^2dxdt\le C.\ea\ee
On the other hand, if $\b\in (0,2),$  multiplying \eqref{1.3}   by $(\te-2)_+\te^{ -\frac{\b}{2}}$ and integration by parts gives
\be \la{p2}\ba  & \left(\int_{\O}\int_2^\te (s-2)_+s^{-\b/2}dsdx \right)_t + \int_{(\te>2)(t)} \frac{\te^{\frac{\b}{2}}\te_x^2}{v}((1-\b/2)+\b\te^{-1})dx\\& =-\int_\O \frac{(\te-2)_+\te^{1-\frac{\b}{2}} }{v}u_xdx+\int_\O \frac{(\te-2)_+\te^{ -\frac{\b}{2}} u_x^2}{v} dx \\& \le C\int_\O (\te-2)_+\te^{2-\frac{\b}{2}}  dx+2\int_\O\frac{(\te-2)_+\te^{ -\frac{\b}{2}} u_x^2}{v} dx \\& \le C\max_{x\in\O}(\te-3/2)^2_+\int_{(\te>2)(t)} \te^{1-\frac{\b}{2}}  dx+C\max_{x\in\O}(\te- 2)_+\int_\O u_x^2  dx\\& \le C\max_{x\in\O}(\te-3/2)^2_+  +C\left(\int_\O u_x^2  dx\right)^2 ,\ea\ee
where in the last inequality we have used \eqref{nep1} and \eqref{nep2}. Since $\b<2,$
it follows from \eqref{p2},  \eqref{vx1}, and \eqref{2.12} that
\bnn \la{sq2}\ba  \int_0^T\int_\O  \te_x^2dxdt& \le  C\int_0^T\int_{(\te<3)(t)} \te^{\b-2}  \te_x^2dxdt+C\int_0^T\int_{(\te>2)(t)}\te^{\b/2}\te_x^2 dxdt \\&
\le C+C \int_0^T\left(\int_{\O} u_x^2dx\right)^2dt ,\ea\enn  which together with   \eqref{le5eq1},    \eqref{w2.11},   \eqref{df8}, and Gronwall's inequality shows that \eqref{lm5w}  still holds.

Finally, it follows from \eqref{mom5}, \eqref{lm5w}, \eqref{w2.11}, \eqref{df8}, \eqref{z4}, and \eqref{vx1} that
\bnn  \int_0^T\int_{\O} u_t^2dxdt\le C,\enn  which together with \eqref{lm5w},  \eqref{w2.11}, and \eqref{df8} gives \eqref{lm5} and finishes the proof of Lemma \ref{lemm5}.\thatsall

Now, we can prove the uniform lower and upper bounds of the temperature $\te.$

\begin{lemma}\la{lmm8} There exists a positive constant $C$ such that for any $(x,t)\in \O\times[0,T]$ \be\ba\la{tq2}  C^{-1}\le \te(x,t)\le C  . \ea\ee\end{lemma}

\pf First, for $p>\b+1,$
multiplying \eqref{1.3} by
$ (\te-2)_+^{p-1}$ and integrating the resultant equality in $x$ over $\O$  leads to
\be\la{p7}\ba&\frac{1}{p}\left(\int_{\O} (\te-2)_+^pdx\right)_t  +(p-1)\int_{\O} \frac{\te^\b(\te-2)_+^{p -2}\te_x^2}{v}dx \\& =\int_{\O} \frac{(\te-2)_+^{p-1}u_x^2}{v}dx-\int_{\O} \frac{\te(\te-2)_+^{p-1}u_x }{v}dx \\& \le C(\ve)\int_{\O} (\te-2)_+^{p-1}u_x^2 dx+\ve \int_\O \te^2(\te-2)_+^{p-1} \\& \le C(\ve)\left(\int_{\O}  u_x^2 dx\right)^{(\b+p+1)/(\b+2)}+C\ve \max_{x\in \O} (\te-3/2)_+^{\b+p+1}  \\& \le C(\ve)\int_{\O}  u_x^2 dx+C\ve \int_\O \te^{\b+p-2}\te_x^2dx  ,\ea\ee where in the last inequality we have used \eqref{lm5} and \eqref{lia22}.

Then, it follows from \eqref{2.12} that
\bnn \ba &\int_0^T\int_\O \te^{\b+p-2}\te_x^2dxdt\\&\le C\int_0^T\int_{\O} \frac{\te^\b(\te-2)_+^{p -2}\te_x^2}{v}dx dt+C\int_0^T\int_{(\te<3)(t)} \te^{\b -2}\te_x^2dxdt\\&\le C\int_0^T\int_{\O} \frac{\te^\b(\te-2)_+^{p -2}\te_x^2}{v}dx dt+C ,\ea \enn which together with
\eqref{p7}, \eqref{w2.11}, and \eqref{df8}    gives \be\ba\la{lm7eq1}\sup_{0\le t\le T}\int_{\O} (\te-2)_+^pdx+\int_0^T\int_{\O} \te^{p+\b-2}\te_x^2dxdt\le C(p) .\ea\ee

Next, multiplying \eqref{1.3} by
$  \te^{\b}\te_t$ and integrating the resultant equality over $\O$ yields
\bnn\ba &  \int_{\O} \te^\b\te_t^2dx+ \int_{\O}\frac{ \te^{\b+1}\theta_tu_x}{v  }dx\\&=\int_{\O}\te^{\b}\te_t\left(\frac{\te^\b\theta_{x}}{v}\right)_{x}dx+\int_{\O}\frac{ \te^{\b }\theta_tu_x^2}{v }dx  \\& =-\int_{\O}\frac{\te^\b\theta_{x}}{v}\left(\te^{\b}\te_t\right)_{x}dx+\int_{\O}\frac{ \te^{\b }\theta_tu_x^2}{v }dx   \\& =-\int_{\O}\frac{\te^\b\theta_{x}}{v}\left(\te^{\b}\te_x\right)_{t}dx+\int_{\O}\frac{ \te^{\b }\theta_tu_x^2}{v }dx  \\&=-\frac{1}{2}
\int_{\O}\frac{\left((\te^\b\theta_{x})^2\right)_t}{v}dx+\int_{\O}\frac{ \te^{\b }\theta_tu_x^2}{v }dx \\&=-\frac{1}{2} \left(\int_{\O}\frac{(\te^\b\theta_{x})^2}{v}dx\right)_t
-\frac{1}{2}\int_{\O}\frac{(\te^\b\theta_{x})^2u_x}{v^2}dx+\int_{\O}\frac{ \te^{\b }\theta_tu_x^2}{v }dx ,\ea\enn
which gives
\be\ba\la{lm8eq2} &  \int_{\O} \te^\b\te_t^2dx+ \frac{1}{2} \left(\int_{\O}\frac{(\te^\b\theta_{x})^2}{v}dx\right)_t\\&=- \frac{1}{2}\int_{\O}\frac{(\te^\b\theta_{x})^2u_x}{v^2}dx-\int_{\O}\frac{ \te^{\b+1}\theta_t u_x}{v  }dx+\int_{\O}\frac{ \te^{\b }\theta_tu_x^2}{v }dx\\&\le C\max_{x\in\O} |u_x|\int_{\O}\left(\te^\b\te_x\right)^2dx +\frac{1}{2}\int_{\O} \te^\b\te_t^2dx +C\int_{\O} \te^{\b+2}u_x^2dx\\&\quad+C\int_{\O}(\te-2)_+^\b u_x^4dx+C\int_{\O}  u_x^4dx\\&\le C\left(\int_{\O}\left(\te^\b\te_x\right)^2dx\right)^2+C \max_{x\in\O} u_x^2+C \max_{x\in\O} u_x^4+\frac{1}{2}\int_{\O} \te^\b\te_t^2dx , \ea\ee
due to \eqref{lm7eq1} and \eqref{lm5}.

Next, it follows from  \eqref{ux2}  and \eqref{lm5}    that
\be\la{tq3} \int_0^T \max_{x\in\O} ( u_x^2+u_x^4)dt\le C,\ee which together with \eqref{lm8eq2},  the Gronwall inequality,  and \eqref{lm7eq1} leads to
\be\ba\la{tebtex}  \sup_{0 \le t\le T}\int_{\O} \left(\te^\b\theta_{x}\right)^2 dx+\int_0^T\int_{\O} \te^\b\te_t^2dxdt\le C. \ea\ee
Combining this with \eqref{nep2} in particular gives
\bnn\ba \max_{x\in\O} (\te-2)_+ & \le\int_{(\te>2)(t)}  |\te_x|dx\\&\le C\left(\int_{(\te>2)(t)}\left(\te^\b\te_x\right)^2dx\right)^{1/2} \le C,\ea\enn which implies that for all $(x,t)\in\O\times [0,\infty),$
\be\la{teup} \te(x,t)\le C.\ee

Next, multiplying  \eqref{1.3} by $(\te-1)^5$ and integrating the resulting equality over $\O$ gives
\be\la{w2.65}\ba &  \frac{1}{6}\left|(\int_\O (\te-1)^6dx)_t\right|\\&=\left|-5\int_\O \frac{(\te-1)^4\te_x^2}{v}dx-\int_\O \frac{\te(\te-1)^{5}u_x}{v}dx +\int_\O \frac{ (\te-1)^{5}u^2_x}{v}dx\right|\\&\le C\int_\O \te_x^2dx+C \int_\O(\te-1)^6dx+C\int_\O u_x^2dx\\&\le C\int_\O \te_x^2dx +C\int_\O u_x^2dx\ea\ee
where in the last inequality we have used \bnn \int_\O (\te-1)^6dx\le C (\int_\O (\te-1)^2dx)^2 \int_\O  \te_x^2dx\le C  \int_\O  \te_x^2dx,\enn
due to \eqref{2.12} and \eqref{teup}. Combining this, \eqref{w2.65}, and \eqref{lm5} gives
\be \la{w2.67}\ba \lim_{t\rightarrow \infty}\int_{\O} \left(\te -1\right)^6dx=0. \ea\ee

Sobolev's inequality shows
\bnn \ba \max_{x\in\O}(\te^{\b+1}-1)^2&\le C\left(\int_\O(\te^{\b+1}-1)^6dx\right)^{1/4} \left(\int_\O\te^{2\b }\te_x^2 dx\right)^{1/4}\\& \le C\left(\int_\O(\te -1)^6dx\right)^{1/4}  , \ea\enn
which together with \eqref{w2.67}    implies that there exists some  $  T_0>0$ such that
\be\la{x5}\ba \te (x,t)\ge 1/2, \ea\ee   for all $(x,t)\in [0,1]\times [ T_0,\infty).$ Moreover, it follows from  \cite[Lemma 2.3]{lsx1} that  there exists some constant $C \ge  2$
such that  \bnn \theta(x,t)\ge  C^{-1} ,\enn  for all $(x,t)\in [0,1]\times [0,T_0].$ Combining this, \eqref{x5}, and \eqref{teup} gives \eqref{tq2}.
The proof of Lemma  \ref{lmm8} is finished. \thatsall

Finally,
we have the following uniform estimate on the $L^2((0,1)\times (0,T))$-norm of $\te_t$ and $\te_{xx}.$

\begin{lemma}\la{tq8}There exists a positive constant $C$ such that \be\la{tq7}\ba \sup_{0\le t\le T}\int_{\O} \te_x^2dx+\int_0^T\int_{\O}\left( \te_t^2+\te_{xx}^2\right)dxdt\le C.\ea\ee\end{lemma}
\pf First, both \eqref{tq2} and \eqref{tebtex} lead  to
\be\ba\la{lm9eq2}  \sup_{0 \le t\le T}\int_{\O} \theta_{x}^2 dx+\int_0^T\int_{\O} \te_t^2dxdt\le C. \ea\ee

Next, it follows from \eqref{1.3} that
\bnn\ba \frac{\left(\te^\b\te_x\right)_x}{v}=  \frac{ \te^\b\te_xv_x}{v^2}- \frac{u_x^2}{v}+ \frac{ \te u_x}{v}+\te_t,\ea\enn
which together with \eqref{z4}, \eqref{tq2},  \eqref{lm5}, \eqref{tq3},   and \eqref{lm9eq2} gives
\be\la{tq4}\ba  \int_0^T\int_{\O}\left| \left(\te^\b\te_x\right)_x\right|^2dxdt  & \le C\int_0^T\max_{x\in\O}\left(  \te^\b\te_x\right)^2\int_{\O} v_x^2dxdt+C\\ & \le C\int_0^T\max_{x\in\O}\left(  \te^\b\te_x\right)^2 dt+C.\ea\ee

We get by \eqref{pj}  and \eqref{tq2}, \bnn\ba \int_0^T\max_{x\in\O}\left(  \te^\b\te_x\right)^2 dt    \le   C   (\de) \int_0^T\int_\O\left(  \te^\b\te_x\right)^2dx  dt +\de \int_0^T\int_{\O}\left| \left(\te^\b\te_x\right)_x\right|^2dxdt ,  \ea\enn  which together with  \eqref{tq4},  \eqref{tq2}, and \eqref{lm5}  implies
\be\la{tq5}\ba
&\int_0^T \max_{x\in\O} \te_x^2dt+\int_0^T\int_{\O}\left| \left(\te^\b\te_x\right)_x\right|^2dxdt  \le C.
\ea\ee

Finally, since
\bnn\te_{xx}=\frac{\left(\te^\b\te_x\right)_x}{\te^\b} -  \frac{\b \te_x^2}{\te } ,\enn it follows from \eqref{tq5},  \eqref{tq2},  and \eqref{tebtex} that
\bnn \ba\int_0^T\int_{\O}  \te_{xx}^2dxdt  &\le  C    \int_0^T\int_{\O}\left| \left(\te^\b\te_x\right)_x\right|^2dxdt+ C \int_0^T\max_x \te_x^2\int_{\O} \te_x^2dxdt  \\&\le C+C\sup_{0\le t\le T}\int_{\O} \te_x^2dx \int_0^T\max_{x\in\O}\te_x^2dt
\\&\le  C, \ea\enn  which together with \eqref{lm9eq2}
gives \eqref{tq7} and finishes the proof of Lemma \ref{tq8}.\thatsall


\end{document}